\documentclass[12pt,twoside,a4paper]{amsart}
\usepackage{amssymb}

\date{\today}

\def\codim{{\rm codim\,}}

\def\dbar{\bar\partial}

\def\C{{\mathbb C}}

\def\D{{\mathcal D}}

\def\L{{\mathcal L}}

\def\Re{{\rm Re\,  }}

\def\L{{\mathcal L}}
\def\U{{\mathcal U}}

\def\be{\begin{equation}}
\def\ee{\end{equation}}

\newtheorem{thm}{Theorem}[section]
\newtheorem{lma}[thm]{Lemma}

\newtheorem{prop}[thm]{Proposition}

\theoremstyle{definition}

\theoremstyle{remark}

\newtheorem{preremark}{Remark}
\newtheorem{preex}{Example}

\newenvironment{remark}{\begin{preremark}}{\qed\end{preremark}}
\newenvironment{ex}{\begin{preex}}{\qed\end{preex}}

\numberwithin{equation}{section}

\begin{document}

\title[Products of residue currents...
]{Products of residue currents of Cauchy-Fantappi\`e-Leray type}

\date{\today}

\author{Elizabeth Wulcan}

\address{Department of Mathematics\\Chalmers University of Technology and the University of G\"oteborg\\S-412 96 G\"OTEBORG\\SWEDEN}

\email{wulcan@math.chalmers.se}



\begin{abstract}
With a given holomorphic section of a Hermitian vector bundle, one can associate a residue current by means of Cauchy-Fantappi\`e-Leray type formulas. In this paper we define products of such residue currents. We prove that, in the case of a complete intersection, the product of the residue currents of a tuple of sections coincides with the residue current of the direct sum of the sections.
\end{abstract}


\maketitle

\section{Introduction}\label{intro}

Let $f$ be a holomorphic function defined in some domain in $\mathbb C^n$ and let $Y=f^{-1}(0)$. Then there exists a distribution $U$ such that $fU=1$, as shown by Schwartz \cite{S}.
For example, one can let $U$ be the principal value distribution $[1/f]$, defined as
\begin{equation*}
\mathcal D_{n,n}\ni\phi \mapsto \lim_{\varepsilon \to 0}\int_{|f|>\varepsilon}\frac{\phi}{f}.
\end{equation*}
The existence of this limit was proven by 
Herrera and Lieberman, \cite{HL}, using Hironaka's desingularization theorem. 
By the Mellin transform, see for example \cite{PT}, one can show that the limit is equal to the analytic continuation to $\lambda=0$ of 
\begin{equation}
\lambda \mapsto \int |f|^{2\lambda} \frac{\phi}{f}.
\end{equation}
The \emph{residue current} associated with $f$ is defined as $\dbar[1/f]$; it has  support on $Y$ and its action on a test form $\phi\in\mathcal D_{n,n-1} $ is given by the analytic continuation to $\lambda=0$ of 
\begin{equation*}
\lambda \mapsto \int\dbar|f|^{2\lambda}\wedge\frac{\phi}{f}.
\end{equation*}

This paper concerns products of residue currents. 
Recall that it is in general not possible to multiply currents (or distributions). However, given a tuple of holomorphic functions $f=(f_1,\ldots, f_m)$, by certain limiting processes one can give meaning to the expression
\begin{equation}\label{coleff}
\dbar\Big[\frac{1}{f_1}\Big ] \wedge\ldots\wedge\dbar\Big [ \frac{1}{f_m}\Big ],
\end{equation}
as was first done by Coleff and Herrera, \cite{CH}. By the Mellin transform, this so called Coleff-Herrera current, denoted by $R^f_{CH}$, can be realized as the analytic continuation to $\lambda=0$ of
\begin{equation}\label{coleffmellin}
\dbar|f_1|^{2\lambda}\frac{1}{f_1}\wedge\ldots\wedge\dbar|f_m|^{2\lambda}\frac{1}{f_m}.
\end{equation}
In case $f$ defines a complete intersection, that is, the codimension of $Y=f^{-1}(0)$ is $m$, then $R^f_{CH}$ has especially nice calculus properties. For example $f_i R^f_{CH}=0$ for all $i$, see \cite{P}, which yields one direction of the duality theorem, due to Passare, \cite{P2}, and Dickenstein-Sessa, \cite{DS}, that asserts that if $f$ is a complete intersection, then a holomorphic function ~$\varphi$ belongs to the ideal $(f)$ if and only if $\varphi R^f_{CH}=0$.

In \cite{PTY} Passare, Tsikh and Yger introduced an alternative approach to multidimensional residue currents by constructing currents based on the Bochner-Martinelli kernel. For each ordered index set $\mathcal I\subseteq \{1,\ldots, m\}$ of cardinality $k$, let $R^f_{\mathcal I}$ be the analytic continuation to $\lambda=0$ of 
\begin{equation*}
\dbar |f|^{2\lambda}\wedge
\sum_{\ell=1}^k(-1)^{k-1}
\frac
{\overline{f_{i_\ell}}\bigwedge_{q\neq \ell} \overline{df_{i_q}}}
{|f|^{2k}},
\end{equation*}
where $|f|^2=|f_1|^2+\ldots+|f_m|^2$.
Then $R^f_\mathcal I$ is a well-defined $(0,k)$-current with support on $Y$, that vanishes whenever $k<\codim Y$ or $k>\min(m,n)$. If $f$ is a complete intersection, there is only one nonvanishing current, namely $R^f_{\{1,\ldots,m\}}$, which corresponds to the classical Bochner-Martinelli kernel and which we denote by $R_{BM}^f$. Then we have the following result.
\begin{thm}[Passare, Tsikh, Yger \cite{PTY}]\label{sondag}
Assume that $f$ is a complete intersection. Then
\[R_{BM}^f=R_{CH}^f.\]
\end{thm}
The Bochner-Martinelli residue currents $R_{\mathcal I}^f$ have been used for investigations in the non-complete intersection case; for example, in \cite{BY}, Berenstein and Yger used them to construct Green currents.

Based on the work in \cite{PTY} Andersson, \cite{A}, introduced more general globally defined residue currents by means of Cauchy-Fantappi\`e-Leray type formulas. Let us briefly recall his construction.
Assume that $f$ is a holomorphic section of the dual bundle $E^*$ of a holomorphic $m$-bundle $E\to X$ over a complex manifold $X$. 
On the exterior algebra over $E$ we have mappings $\delta_f: \Lambda^{\ell+1} E \to \Lambda^\ell E$ of interior multiplication by ~$f$, and $\delta_f^2=0$. 
Let $\mathcal E_{0,k}(X,\Lambda^\ell E)$ be the space of smooth sections of the exterior algebra of $E^*\oplus T_{0,1}^*$ which are $(0,k)$-forms with values in $\Lambda^\ell E$, and let $\mathcal D'_{0,k}(X,\Lambda^\ell E)$ be the corresponding space of currents. The mappings $\delta_f$ extend to these spaces, where it anti-commutes with ~$\dbar$. Thus $\mathcal D'_{0,k}(X,\Lambda^\ell E)$ is a double complex and the corresponding total complex is
\begin{equation*}
\cdots \stackrel{\nabla_f}{\longrightarrow} \L^{r-1}(X,E)
\stackrel{\nabla_f}{\longrightarrow}\L^{r}(X,E)
\stackrel{\nabla_f}{\longrightarrow}\cdots,
\end{equation*}
where $\L^r(X,E)=\bigoplus_{k+\ell=r}\D'_{0,k}(X,\Lambda^{-\ell} E)$ and $\nabla_f=\delta_f-\dbar$. The exterior product, $\wedge$, induces a mapping
\begin{equation*}
\wedge: \L^r(X,E)\times \L^s(X,E)\to \L^{r+s}(X,E)
\end{equation*}
when possible, and $\nabla_f$ is an antiderivation with respect to $\wedge$.

If $\varphi$ is a holomorphic function such that is $\varphi=\nabla_f v$ for some $v\in \L^{-1}(X,E)$, one can prove, provided $X$ is Stein, that there is a holomorphic solution $\psi$ to the division problem $\sum\psi_j f_j=\varphi$.
Andersson's idea to find such a $v$ was to start looking for a solution to $\nabla_f u=1$. Assume that $E$ is equipped with some Hermitian metric and let $s$ be the section of ~$E$ with pointwise minimal norm such that $\delta_f s=|f|^2$ and let 
\begin{equation*}
u^f =
\frac{s}{\nabla_{f} s}= \frac{s}{\delta_{f} s-\dbar s}=
\sum_\ell\frac{s\wedge(\dbar s)^{\ell-1}}{(\delta_f s)^{\ell}}=
\sum_\ell\frac{s\wedge(\dbar s)^{\ell-1}}{|f|^{2\ell}}
\end{equation*}
be the Cauchy-Fantappi\`e-Leray form, introduced in \cite{A2} in order to construct integral formulas in a convenient way. Clearly $u^f\in \L^{-1}$ is well-defined outside $Y$ and since $\nabla_f s$ is of even degree the expression $s/\nabla_{f} s$ makes sense, and it follows that $\nabla_f u^f=1$ outside $Y$.  
In \cite{A} it is proved that the form $|f|^{2\lambda} u^f$ has an analytic continuation as a current to $\Re\lambda>-\epsilon$. The value at $\lambda=0$, denoted by ~$U^f$, yields an extension of $u^f$ over $Y$. In analogy with the one function case, we will sometimes refer to $U^f$ as the principal value current. Clearly, if $Y\neq\emptyset$, $U^f$ can not fulfill $\nabla_f U^f$. In fact, $\nabla_f U^f=1-R^f$, where $R^f=\dbar|f|^{2\lambda}\wedge u^f|_{\lambda=0}$ now defines the residue current of $f$. It holds that $R^f=R_p+\ldots+R_\mu$, where $R_j\in \mathcal D'_{0,j}(X,\Lambda^j)$, $p=\codim Y$ and $\mu=\min(m,n)$. 
Moreover, if $\varphi R^f=0$, then $v=Uf$ yields the desired solution to $\nabla_f v=\varphi$ and thus ~$\varphi$ belongs to the ideal generated by $f$ locally.

If $E$ is a trivial bundle endowed with the trivial metric, the coefficients of $R^f$ will actually be the Bochner-Martinelli currents $R^f_{\mathcal I}$. If $f$ is a complete intersection, the only nonvanishing coefficient will be $R_{BM}^f$.

Our first goal is to define products of currents of the type $U^f$ and ~$R^f$.
Let us consider \eqref{coleffmellin}. If we assume that each
$f_i$ is a section of the dual bundle $E_i^*$ of a line bundle $E_i$ with frame $e_i$ and dual frame $e_i^*$, the Cauchy-Fantappi\`e-Leray form $u^{f_i}$ is just $e_i/f_i$, so in fact \eqref{coleffmellin} times the element $e_1\wedge\ldots\wedge e_r$ can be expressed as
\begin{equation}\label{produkten}
\dbar|f_1|^{2\lambda} \wedge u^{f_1}\wedge\ldots\wedge
\dbar|f_r|^{2\lambda} \wedge u^{f_r}.
\end{equation}
In light of this, it is most tempting to extend this product to include not only sections of line bundles but sections $f_i$ of bundles of arbitrary rank. 
To be more accurate, we assume that $f_i$ is a section of the dual bundle of a holomorphic $m_i$-bundle $E_i\to X$. 
Further, we assume that each $E_i$ is equipped with a Hermitian metric, we let $s_i$ be the section of $E_i$ of minimal norm such that $\delta_{f_i} s_i=|f_i|^2$, and we let $u^{f_i}$ be the corresponding Cauchy-Fantappi\`e-Leray form. Then 
\eqref{produkten} has meaning as a form taking values in the exterior algebra over $E=E_1\oplus \cdots \oplus E_r$. 
Thus, in accordance with the line bundle case, we can take the value at $\lambda=0$ of \eqref{produkten} as a definition of $R^{f_1}\wedge\ldots\wedge R^{f_r}$, provided that the analytic continuation exists. However, this is assured by Proposition ~\ref{existensen}, where products are defined also of
principal value currents.

\begin{prop}\label{existensen}
Let $f_i$ be holomorphic sections of the Hermitian 
\linebreak
$m_i$-bundles $E_i^*\to X$. Let $u^{f_i}$ be the corresponding Cauchy-Fantappi\`e-Leray forms and let $Y_i=f_i^{-1}(0)$. 
Then 
\begin{equation}\label{form_existensen}
\lambda \mapsto
|f_r|^{2\lambda}u^{f_r} \wedge 
\ldots \wedge
|f_{s+1}|^{2\lambda}u^{f_{s+1}} \wedge 
\dbar|f_s|^{2\lambda} \wedge u^{f_s} \wedge 
\ldots \wedge
\dbar|f_1|^{2\lambda}\wedge u^{f_1}
\end{equation}
has an analytic continuation as a current to $\lambda > - \epsilon$. 

We define $T=U^{f_r}\wedge \ldots \wedge U^{f_{s+1}} \wedge R^{f_s} \wedge \ldots \wedge R^{f_1}$ as the value at $\lambda=0$. Then $T$ has support on $\bigcap_{i=1}^s Y_i$ and it is alternating with respect to the principal value factors $U^{f_i}$ and symmetric with respect to the residue factors $R^{f_i}$.
\end{prop}

Of course there is nothing special about the ordering that we have chosen; we can just as well mix $U$'s and $R$'s. 

If the bundle $E$ is trivial, endowed with the trivial metric, and moreover if $f_1\oplus\cdots\oplus f_r$ is a complete intersection, then $R^{f_r}\wedge \ldots\wedge R^{f_1}$ will consist of only one term, which can be interpreted as a product of the corresponding Bochner-Martinelli currents $R_{BM}^{f_i}$. 
In general, however, there will also occur terms of lower degree. 
\begin{prop}\label{tove}
Let 
\[T=U^{f_r} \wedge \ldots \wedge U^{f_{s+1}} \wedge R^{f_s} \wedge \ldots \wedge R^{f_1}\]
be defined as above. Let $m=m_1+\ldots + m_r$. Then $T=T_{p} + \ldots + T_{q}$, where $T_\ell \in \D'_{0,\ell}(\Lambda^\bullet E)$, $p=\codim Y_1 \cap \ldots \cap Y_s$ and $q=\min(m,n)$.
In particular, if $f=f_1\oplus\cdots\oplus f_r$ is a complete intersection, then
$R^{f_r} \wedge \ldots \wedge R^{f_1}$
consists of only one term of top degree $m$.
\end{prop}

Our next aim is to prove a generalized version of Theorem \ref{sondag}. 
Since, in the particular case when the bundles $E_i$ are all line bundles, the current
$R^{f_1}\wedge\ldots\wedge R^{f_r}$
is just the Coleff-Herrera current of $f$ times $e_1\wedge\ldots\wedge e_r$, we can formulate the equivalence in the theorem as
\begin{equation}\label{sondag-bundle}
R^{f_1\oplus\ldots\oplus f_r}=R^{f_1}\wedge\ldots\wedge R^{f_r}.
\end{equation}
Now, the obvious question is, does this equality extend to hold for sections of vector bundles of arbitrary rank.
Our main result states that this is indeed the case.

\begin{thm}\label{main_thm}
Let $f_i$ be holomorphic sections of the Hermitian $m_i$-bundles $E_i^*$ and let $f$ be the section $f_1\oplus\cdots\oplus f_r$ of $E^*=E_1^*\oplus\cdots\oplus E_r^*$. If $f$
is a complete intersection, that is, $\codim f^{-1}(0)=m_1+\ldots+m_r$, then
\[R^f=R^{f_1}\wedge\ldots\wedge R^{f_r}.\]
\end{thm}

That is, in a local perspective, given a tuple of functions split into subtuples, the product of the Bochner-Martinelli 
currents of each subtuple is equal to the Bochner-Martinelli current of the whole tuple of functions.
We give an explicit proof of Theorem ~\ref{main_thm} based on the existence of two $\nabla_f$-potentials.
\begin{thm}\label{potential}
Let $f=f_1\oplus\cdots\oplus f_r$ be a section of $E^*=E_1^*\oplus\cdots\oplus E_r^*$. Assume that $f$ is a complete intersection. Then there exists a current ~$V$ such that 
\begin{equation}\label{potential1}
\nabla_f V = 1-R^{f_1}\wedge\ldots\wedge R^{f_r},
\end{equation}
and furthermore a current $U^f\wedge V$ such that
\begin{equation}
\nabla_f(U^f\wedge V)=V-U^f.
\end{equation}
\end{thm}
At first it might seem a bit peculiar to denote the second potential by $U^f\wedge V$. However, notice that on a formal level, if we were allowed to multiply currents so that $\nabla_f$ acted as an antiderivation on the products, then
\[\nabla_f(U^f\wedge V)=
(1-R^f)\wedge V - U^f\wedge(1-R^{f_1}\wedge\ldots\wedge R^{f_r}),\]
since $U^f$ is of odd degree. 
From Proposition ~\ref{tove} we know that $R^f$ and 
\linebreak
$R^{f_1}\wedge\ldots\wedge R^{f_r}$
take values in $\Lambda^m E$, since $f$ is a complete intersection. But since $V$ and $U^f$ have positive degree in $e_j$ it is reasonable to expect the products $V\wedge R^f$ and $U^f\wedge R^{f_1}\wedge\ldots\wedge R^{f_r}$ to vanish.
Thus we are left with $V-U^f$, and the notation is motivated.

\begin{proof}[Proof of Theorem ~\ref{main_thm}]
Recall that $\nabla_f U^f=1-R^f$.
Hence, applying ~$\nabla_f$ twice to $U^f\wedge V$ yields
\[ 
0=\nabla_f^2(U^f\wedge V)=\nabla_f(U^f-V)=R^{f_1}\wedge\ldots\wedge R^{f_r}-R^f,\]
and thus we are done.
\end{proof}

The disposition of this paper is as follows. In Section ~\ref{CFL-currents} we give proofs of Proposition ~\ref{existensen} and Proposition ~\ref{tove}. In Section ~\ref{complete_intersection} we prove Theorem ~\ref{potential}.
Finally, in Section ~\ref{examples} we give an example of products of Cauchy-Fantappi\`e-Leray currents and also discuss a possible generalization of Theorem \ref{main_thm}. 

\section{Products of residue currents of \\ Cauchy-Fantappi\`e-Leray type}\label{CFL-currents}
We start with the proof of Proposition ~\ref{existensen}.
For further use a slightly more general formulation is appropriate. Indeed, the proof of Theorem ~\ref{potential} requires a  broader definition of products of currents. We need to allow also products of currents of sections of the bundle $E$, that are not necessarily orthogonal, at least in certain cases.
Thus we give a new, somewhat unwieldy, version of Proposition ~\ref{existensen} that however covers all the currents that we will be concerned with.  

By the notion that a form (or current) is of degree $k$ in $d\bar z_j$, we will just mean that it is a 
$(\bullet, k)$-form. In the same manner, we will say that a form is of degree $\ell$ in $e_j$ when it takes values in $\Lambda^\ell E$.

\begin{prop}\label{existensen2}
Let $f=f_1\oplus\ldots\oplus f_r$ be a holomorphic section of the bundle
$E^*=E_1^*\oplus\ldots\oplus E_r^*$, where $E_i^*$ is a Hermitian $m_i$-bundle.
For a subset $I=\{I_1,\ldots,I_p\}$ of $\{1,\ldots,r\}$, 
let $f_I$ denote the section $f_{I_1}\oplus\ldots\oplus f_{I_p}$ of $E_I^*=E_{I_1}^*\oplus\ldots\oplus E_{I_p}^*$, 
let $u^{f_{I}}$ be the corresponding Cauchy-Fantappi\`e-Leray form,
let $Y_{I}=f_{I}^{-1}(0)$,
and let $m_I=m_{I_1}+\ldots + m_{I_p}$. 
If $I^1, \ldots, I^t$ are subsets of $\{1,\ldots,r\}$, then 
\begin{equation}\label{form_existensen2}
\lambda \mapsto
|f_{I^t}|^{2\lambda}u^{f_{I^t}} \wedge 
\ldots \wedge
|f_{I^{s+1}}|^{2\lambda}u^{f_{I^{s+1}}} \wedge 
\dbar|f_{I^s}|^{2\lambda} \wedge u^{f_{I^s}} \wedge 
\ldots \wedge
\dbar|f_{I^1}|^{2\lambda}\wedge u^{f_{I^1}}
\end{equation}
has an analytic continuation to $\lambda > - \epsilon$.
 
We define $T=U^{f_{I^t}}\wedge \ldots \wedge U^{f_{I^{s+1}}} \wedge R^{f_{I^s}} \wedge \ldots \wedge R^{f_{I^1}}$ as the value at $\lambda=0$. Then $T$ has support on $\bigcap_{i=1}^s Y_{I^k}$ and it is alternating with respect to the principal value factors $U$ and commutative with respect to the residue factors $R$.
\end{prop}
Note that Proposition ~\ref{existensen} corresponds to the particular case when each ~$I^j$ is just a singleton.
The proof of Proposition ~\ref{existensen2} is very much inspired by the proof of Lemma 2.2 in \cite{PTY} and Theorem 1.1 in \cite{A}. It is based on the possibility of resolving singularities by Hironaka's theorem, see \cite{At}, and the following lemma, which is proven essentially by integration by parts.
\begin{lma}\label{main_lemma}
Let $v$ be a strictly positive smooth function in $\C$,  
 $\varphi$ a test function in $\C$, and 
$p$  a positive integer.  Then
\[
\lambda\mapsto\int v^\lambda |s|^{2\lambda}\varphi(s)\frac{ds\wedge d\bar s}{s^p}
\]
and
\[
\lambda\mapsto\int \dbar (v^\lambda |s|^{2\lambda})
\wedge \varphi(s)\frac{ds}{s^p}
\]
both have meromorphic  continuations  to the entire plane with poles
at rational points on the negative real axis. 
At $\lambda=0$ they are both independent of $v$, and the
second one only depends on the germ of $\varphi$ at the origin.
Moreover, if $\varphi(s)=\bar s\psi(s)$ or $\varphi=d\bar s\wedge\psi$, then
the value of the second integral at $\lambda=0$ is zero.
\end{lma}

\begin{proof}[Proof of Proposition ~\ref{existensen2}]
We may assume that the bundle $E=E_1\oplus\cdots\oplus E_r$ is trivial since the statement is clearly local. Note that $f_i=\sum f_{i,j}e_{i,j}^*$, where $e_{i,j}^*$ is the trivial frame.
The proof is based on the possibility to resolve singularities locally using Hironaka's theorem. Given a small enough neighborhood $\U$ of a given point in $X$ there exist a $n$-dimensional manifold ~$\widetilde \U$ and a proper analytic map $\Pi_h:\widetilde \U \to \U$ such that if $Z=\{\prod_{i,j} f_{i,j}=0\}$ and $\widetilde Z =\Pi_h^{-1}(Z)$, then $\Pi:\widetilde \U \setminus \widetilde Z \to \U \setminus Z$ is biholomorphic and such that moreover $\widetilde Z$ has normal crossings in $\widetilde \U$. 
This implies that locally in $\widetilde \U$ we have that $\Pi_h^* f_{i,j}=a_{i,j}\mu_{i,j}$, where ~$a_{i,j}$ are non-vanishing and $\mu_{i,j}$ are monomials in some local coordinates $\tau_k$. 
Further, given a finite number of monomials $\mu_1\ldots,\mu_m$ in some coordinates $\tau_k$ defined in an $n$-dimensional manifold $\U_t$, there exists a toric manifold $\widetilde \U_t$ and a proper analytic map $\Pi_t\colon \widetilde \U_t\to \U_t$ such that $\Pi_t$ is biholomorphic outside the coordinate axes and moreover, locally it holds that, for some $i$, $\Pi_t^*\mu_i$ divides all $\Pi_t^*\mu_j$, see \cite{BGVY} and \cite{KKMS}. Clearly, if $\mu_i$ divides $\mu_j$ in $\U_t$ then $\Pi_t^*\mu_i$ divides $\Pi_t^*\mu_j$ in $\widetilde \U_t$. Thus after a number, say ~$q$, of such toric resolutions $\Pi_{t_i}$ we can locally consider each section $f_{I^i}$ as a monomial times a non-vanishing section. More precisely we have that $\Pi^*f_{I^i}=\mu_i f'_{I^i}$, where $\Pi=\Pi_{t_q}\circ\cdots\circ\Pi_{t_1}\circ\Pi_h$, $\mu_i$ is a monomial and $f'_{I^i}$ is a non-vanishing section of $E_{I^i}^*$.

Let $\phi$ be a test form with compact support. After a partition of unity we may assume that it has support in a neighborhood $\U$ as above. Then, since $\Pi_h$ is proper, the support of $\Pi_h^*\phi$ can be covered by a finite number of neighborhoods in which it holds that $\Pi_h^*\phi=a_{i,j}\mu_{i,j}$. If $\psi$ is a test form with support in such a neighborhood, then the support of $\Pi_{t_1}^*\psi$ can be covered by finitely many neighborhoods in which we have the desired property that the pull-back of one monomial divides some of the other ones, and so on. Thus, for $\Re \lambda> 2 \max_i m_{I^i}$, \eqref{form_existensen2} is in $L^1_\text{loc}$, and 
since $\Pi$ is biholomorphic outside a set of measure zero we have that 
\begin{equation*}
\int 
|f_{I^t}|^{2\lambda}u^{f_{I^t}} \wedge 
\ldots \wedge
|f_{I^{s+1}}|^{2\lambda}u^{f_{I^{s+1}}} \wedge 
\dbar|f_{I^s}|^{2\lambda}\wedge u^{f_{I^s}} \wedge 
\ldots \wedge
\dbar|f_{I^1}|^{2\lambda}\wedge u^{f_{I^1}}
\wedge \phi 
\end{equation*}
is equal to a finite number of integrals of the form
\begin{equation}\label{mandag}
\int \Pi^*(|f_{I^t}|^{2\lambda}u^{f_{I^t}} \wedge 
\ldots \wedge
|f_{I^{s+1}}|^{2\lambda}u^{f_{I^{s+1}}} \wedge 
\dbar|f_{I^s}|^{2\lambda}\wedge u^{f_{I^s}} \wedge 
\ldots \wedge
\dbar|f_{I^1}|^{2\lambda}\wedge u^{f_{I^1}})
\wedge ~\tilde \phi.
\end{equation}
Here 
\[
\tilde\phi=
\rho_{t_q}\Pi_{t_q}^*(\ldots
\rho_{t_{1}}\Pi_{t_{1}}^*(\rho_h \Pi_h^*(\phi))),
\]
where the $\rho_\bullet$'s are functions from some partitions of unity, so that the test form $\tilde\phi$ has support in a neighborhood where it holds that $\Pi^* f_{I^i}= \mu_i f_{I^i}'$. In such a coordinate neighborhood the pullback of ~$s_{I^i}$ is ~$\bar \mu_i$ times a smooth form, so that 
$\Pi^*(s_{I^i} \wedge (\dbar s_{I^i})^{\ell-1})$ is $\bar \mu_i^\ell$ times a smooth form. Moreover $\Pi^*|f_{I^i}|^2 = |\mu_i|^2 a_i$, where $a_i$ is a strictly positive smooth function.
Thus 
\[
\Pi^* u^{f_{I^i}}=\sum_\ell \frac{\bar \mu_i^\ell \alpha_{i,\ell}}{|\mu_i|^{2\ell}}
=\sum_\ell \frac{\alpha_{i,\ell}}{\mu_i^\ell},
\] 
where ${\alpha_{i,\ell}}$ are smooth forms taking values in $\Lambda^\ell E$, and so \eqref{mandag} is equal to a finite sum of integrals
\begin{multline}\label{int_erik}
\int
|\mu_t|^{2\lambda} a_t^{\lambda} \frac{\alpha_{t,\ell_t}}{\mu_t^{\ell_t}} \wedge 
\ldots \wedge
|\mu_{s+1}|^{2\lambda} a_{s+1}^{\lambda} \frac{\alpha_{s+1,\ell_{s+1}}}{\mu_{s+1}^{\ell_{s+1}}} \wedge
\\
\dbar (|\mu_s|^{2\lambda} a_s^{\lambda}) \wedge \frac{\alpha_{s,\ell_s}}{\mu_s^{\ell_s}} \wedge 
\ldots \wedge
\dbar (|\mu_1|^{2\lambda} a_1^{\lambda}) \wedge \frac{\alpha_{1,\ell_1}}{\mu_1^{\ell_1}} \wedge
\tilde \phi.
\end{multline}
Expanding each factor $\dbar(|\mu_j|^{2\lambda}a_j^\lambda)$ by Leibniz' rule results in a finite sum of terms. Letting $\dbar$ fall only on the monomials $\mu_i$ yields integrals of the form
\begin{equation}\label{int_emil}
\int
a^\lambda |\mu'|^{2\lambda} 
\frac{\alpha_L}{\mu_L} \wedge
\dbar|\sigma_s^{q_s}|^{2\lambda} \wedge
\ldots \wedge
\dbar|\sigma_1^{q_1}|^{2\lambda} \wedge
\tilde \phi,
\end{equation}
where $\sigma_i$ is one of the coordinate functions $\tau_j$ that divide $\mu_i$, $a = a_t \cdots a_1$ is a strictly positive smooth function, $\mu'$ and $\mu_L'$ are monomials in $\tau_j$ not divisible by any $\sigma_i$ and $\alpha_L=C \alpha_{t,\ell_t} \wedge \ldots \wedge \alpha_{1,\ell_1}$ is a smooth form, where $C$ is just a constant that depends on the relation between $q_i$ and the number of ~$\sigma_i$'s in $\mu_i$.
The remaining integrals, that arise when $\dbar$ falls on any of the $a_i$, 
vanish in accordance with Lemma ~\ref{main_lemma}. Indeed, consider one of the integrals obtained when $\dbar$ falls on $a_1$,
\begin{equation*}
\lambda ~~ 
\int
a^\lambda |\mu'|^{2\lambda} 
\frac{\alpha_L}{\mu_L} \wedge
\dbar|\sigma_s^{q_s}|^{2\lambda} \wedge
\ldots \wedge
\dbar|\sigma_2^{q_2}|^{2\lambda} \wedge
\dbar a_1 \wedge
\tilde \phi.
\end{equation*}
This is just $\lambda$ times an integral of the form \eqref{int_emil}, so provided that we can prove the existence of an analytic continuation of \eqref{int_emil}, it must clearly vanish at $\lambda=0$. 
 
Now an application of Lemma ~\ref{main_lemma} for each $\tau_k$ that divides any of the ~$\mu_j$'s gives the desired analytic continuation of \eqref{int_emil} to $\lambda > -\epsilon$. 
Note that for $\sigma_1, \ldots, \sigma_s$ we get integrals of the second type, for the remaining $\tau_i$ integrals of the first type, so that the value at $\lambda=0$ is a current with support on
$\{\sigma_s=0\} \cap \ldots \cap \{\sigma_1=0\}.$
Thus the value of \eqref{int_erik} at $\lambda=0$ has support on 
\[\{\mu_s=0\} \cap \ldots \cap \{\mu_1=0\}=\widetilde Y_{I^s} \cap \ldots \cap \widetilde Y_{I^1},\]
where $\widetilde Y_\bullet = \Pi^{-1} Y_\bullet$,  
and accordingly $U^{f_{I^t}}\wedge \ldots \wedge U^{f_{I^{s+1}}}\wedge R^{f_{I^s}} \wedge \ldots \wedge R^{f_{I^1}}$ is a current with support on $Y_{I^s}\cap \ldots \cap Y_{I^1}$.

Since the form \eqref{form_existensen2} is alternating with respect to the factors $|f_{I^i}|^{2\lambda} u^{f_{I^i}}$ and symmetric with respect to the factors $\dbar|f_{I^i}|^{2\lambda} \wedge u^{f_{I^i}}$, it follows that $U^{f_{I^t}}\wedge \ldots \wedge U^{f_{I^{s+1}}} \wedge R^{f_{I^s}} \wedge \ldots \wedge R^{f_{I^1}}$ is alternating with respect to the principal value factors and symmetric with respect to the residue factors.
\end{proof}

We continue with the proof of Proposition ~\ref{tove}.

\begin{proof}[Proof of Proposition ~\ref{tove}]
Notice that $T_{\ell}$ is the analytic continuation to 
\linebreak
$\lambda=0$ of the terms 
\begin{equation}\label{form_tove}
|f_r|^{2\lambda}u^{f_r}_{\ell_r} \wedge
\ldots \wedge
|f_{s+1}|^{2\lambda}u^{f_{s+1}}_{\ell_{s+1}} \wedge 
\dbar|f_s|^{2\lambda}\wedge
u^{f_s}_{\ell_s} \wedge 
\ldots \wedge
\dbar |f_1|^{2\lambda}\wedge u^{f_1}_{\ell_1},
\end{equation}
where 
\[u^{f_i}_{\ell_i} = \frac{s_i \wedge (\dbar s_i)^{\ell_i -1}}{|f_i|^{2\ell_i}}\]
and the total degree in $d\bar z_j$ (that is $\ell_1+\ldots +\ell_r-r+s$) is $\ell$.

Following the proof of Proposition ~\ref{existensen2}, a term of the form \eqref{form_tove}, integrated against a test form $\phi$, is equal to a sum of terms like
\begin{multline}\label{int_tove}
\int
|\mu_r|^{2\lambda} a_r^{\lambda} \frac{\alpha_{r,\ell_r}}{\mu_r^{\ell_r}} \wedge 
\ldots \wedge
|\mu_{s+1}|^{2\lambda} a_{s+1}^{\lambda} \frac{\alpha_{s+1,\ell_{s+1}}}{\mu_{s+1}^{\ell_{s+1}}} \wedge
\\
\dbar (|\mu_s|^{2\lambda} a_s^{\lambda}) \wedge \frac{\alpha_{s,\ell_s}}{\mu_s^{\ell_s}} \wedge 
\ldots \wedge
\dbar (|\mu_1|^{2\lambda} a_1^{\lambda}) \wedge \frac{\alpha_{1,\ell_1}}{\mu_1^{\ell_1}} \wedge
\tilde \phi,
\end{multline}
where the $\alpha_{i,\ell_i}$'s are smooth forms of degree $\ell_i$ in $e_j$, the $a_i$'s are non-vanishing functions, the $\mu_i$'s are monomials in some local coordinates ~$\tau_j$ and $\tilde\phi$ is as in the previous proof.
We can find a toric resolution such that locally one of $\mu_1, \ldots, \mu_s$ divides the other ones, so without loss of generality we may assume that $\mu_1$ divides $\mu_2, \ldots, \mu_s$. 

We expand $\dbar(|\mu_1|^{2\lambda}a_1^\lambda)$ by Leibniz' rule. Observe that when $\dbar$ falls on ~$a_1^\lambda$ the integral vanishes as in the proof of Proposition ~\ref{existensen2}, and thus it suffices to consider the case when $\dbar$ falls on one of the $\tau_j$ that divide ~$\mu_1$, say on ~$|\sigma|^{2\lambda}$. If $\ell<p$, we claim that this part of \eqref{int_tove} vanishes when integrating with respect to $\sigma$. In fact, we may assume that 
$\phi=\phi_I\wedge d\bar z_I$, 
where $\phi_I$ is an $(n,0)$-form and 
$d\bar z_I=d\bar z_{I_1}\wedge\ldots\wedge d\bar z_{I_{n-\ell}}$. 
Now $ d\bar z_I$ vanishes on the variety $Y_1 \cap \ldots \cap Y_s$ of codimension $p$ for degree reasons. Consequently $\Pi^*(d\bar z_I)$ vanishes on $\widetilde Y_1 \cap \ldots \cap \widetilde Y_s$, and in particular on 
\linebreak
$\{\sigma=0\}$. However, this is a form in $d\bar\tau_k$ with antiholomorphic
coefficients since $\Pi$ is holomorphic, and therefore each of its terms contains a factor
$d\bar\sigma$ or a factor $\bar\sigma$.
Indeed, if $\Psi(\tau)$ is a form in $d\bar\tau_k$ with antiholomorphic coefficients we can write
\[\Psi(\tau)=\Psi'(\tau)\wedge d\bar\sigma+\Psi''(\tau),\]
where $\Psi''(\tau)$ does not contain $d\bar\sigma$. The first term clearly vanishes on $\{\sigma=0\}$ since $d\bar\sigma$ does. If $\Psi(\tau)$ vanishes on $\{\sigma=0\}$, then $\Psi''(\tau)$ does, and hence it contains a factor $\bar\sigma$ due to antiholomorphicity. In both cases the $\sigma$-integral, and thereby \eqref{int_tove}, vanishes according to Lemma ~\ref{main_lemma}.
\end{proof}

\section{The complete intersection case}\label{complete_intersection}
If $f=f_1,f_2$ defines a complete intersection, then 
\[
\lambda=(\lambda_1,\lambda_2)\mapsto\int\dbar |f_1|^{2\lambda_1}\frac{1}{f_1}\wedge
\dbar |f_2|^{2\lambda_2}\frac{1}{f_2}\wedge \phi
\]
is holomorphic at $\lambda=0$, see \cite{PT}, as was first proven by Berenstein and Yger. Although not yet satisfactorily proven the result is believed to extend to any finite number of functions $f_i$. 

Thus, taking a closer look at the definition of the currents in Proposition ~\ref{existensen} (and Proposition ~\ref{existensen2}) a natural question is whether 
\begin{equation}\label{nagot}
t(\lambda):=
\dbar|f_r|^{2\lambda_r}\wedge u^{f_r} \wedge 
\ldots \wedge
\dbar|f_1|^{2\lambda_1}\wedge u^{f_1}
\end{equation}
is holomorphic in $\lambda$.

In general, letting $\lambda$ tend to 0 along different paths yields different currents as shown by the following example. 
Let $f_1=z_1$ and $f_2=z_1 z_2$ in $\C^2$ and consider $t(\lambda)$ acting on a test form $\phi=\varphi(z_1,z_2) ~dz_1\wedge dz_2$
\begin{equation*}
\int
\frac{\dbar|f_1|^{2\lambda_1}\wedge \dbar|f_2|^{2\lambda_2}}{f_1f_2}\wedge \phi
=
\int
\frac{\dbar|z_1|^{2\lambda_1}\wedge \dbar|z_1 z_2|^{2 \lambda_2}}
{z_1^2 z_2}\wedge \phi.
\end{equation*}
As $\lambda$ tends to zero the integral approaches the value 
$\lambda_1/(\lambda_1+\lambda_2)\varphi_{z_1}(0,0)$. Thus the value at 0 depends on the ratio between $\lambda_1$ and $\lambda_2$. 

In case we have two functions defining a complete intersection though, this phenomenon does not occur. By an integration by parts, we can write ~$t(\lambda)$ integrated against a test form $\phi$ as 
\[
\int
|f_2|^{2\lambda_2} u^{f_2}\wedge 
\dbar|f_1|^{2\lambda_1}\wedge u^{f_1}\wedge \dbar \phi.
\]
After a resolution of singularities this is equal to a sum of integrals of the form
\[
\int
|\mu_2|^{2\lambda_2} a_2^{\lambda_2} 
\frac{\alpha_{2,\ell_2}}{\mu_2^{\ell_2}}\wedge
\dbar(|\mu_1|^{2\lambda_1} a_1^{\lambda_1})\wedge 
\frac{\alpha_{1,\ell_1}}{\mu_1^{\ell_1}}\wedge
\dbar ~\tilde\phi,
\]
according to the proof of Proposition ~\ref{existensen2}.
Let $\sigma$ be one of the coordinate functions that divides $\mu_1$. If $\mu_2$ does not contain $\sigma$, then the $\sigma$-integral is clearly independent of $\lambda$. If $\mu_2$ does contain $\sigma$, at first we seem to end up in a situation similar to the one above, where the result at $\lambda=0$ depends on the relation between $\lambda_1$ and $\lambda_2$. However, the $\sigma$-integral vanishes for the same reason as $T_\ell$ vanishes when $\ell<p$ in the proof of Proposition ~\ref{tove}. Thus, in this case, the definition of the current $T=t(\lambda)|_{\lambda=0}$ is robust in the sense that it does not depend on the particular path along which $\lambda$ tends to zero.

Provided that the Mellin transform of the residue integral is holomorphic in $\lambda$ in a neighborhood of $0 \in \C^r$, it is reasonable to believe that also $t(\lambda)$ is. Presuming this to be true, we can give a soft proof of Theorem ~\ref{main_thm} based on Theorem ~\ref{sondag}.
Indeed, if $f=\sum_1^m f_j e_j^*$ is a section of a bundle $E^*$ we 
let 
\[
t_{CFL}^f(\lambda)=\dbar|f|^{2\lambda} \wedge u^f,
\]
and
\[
t_{CH}^f(\lambda)=\dbar|f_1|^{2\lambda}\frac{1}{f_1}\wedge \ldots \wedge
\dbar|f_m|^{2\lambda}\frac{1}{f_m},
\]
where $CFL$ and $CH$ of course stand for Cauchy-Fantappi\`e-Leray and Coleff-Herrera, respectively.
With this notation the equality in Theorem \ref{sondag} can be expressed as
\begin{equation}\label{penna}
t_{CFL}^f(\lambda)|_{\lambda=0}=t_{CH}^f(\lambda)|_{\lambda=0}.
\end{equation}

Now let $f$ and $g$ be sections of the bundles $E_1^*$ and $E_2^*$, respectively, and assume that $f\oplus g$ is a complete intersection. 
By definition,
\[
R^f\wedge R^g=t_{CFL}^f(\lambda)\wedge t_{CFL}^g(\lambda)|_{\lambda=0},
\]
and 
\[
R^{f\oplus g}=t_{CFL}^{f\oplus g}(\lambda)|_{\lambda=0},
\]
so we need to prove that
\[
t_{CFL}^f(\lambda)\wedge t_{CFL}^g(\lambda)|_{\lambda=0}=
t_{CFL}^{f\oplus g}(\lambda)|_{\lambda=0}.
\]
If $\Re \lambda_2$ is large enough, $t_{CFL}^g(\lambda_2)$ is in $\L^1_\text{loc}$, and so by \eqref{penna}
\[
t_{CFL}^f(\lambda_1)\wedge t_{CFL}^g(\lambda_2)|_{\lambda_1=0}=
t_{CH}^f(\lambda_1)\wedge t_{CFL}^g(\lambda_2)|_{\lambda_1=0},
\]
and analogously, if $\Re \lambda_1$ is large enough
\[
t_{CH}^f(\lambda_1)\wedge t_{CFL}^g(\lambda_2)|_{\lambda_2=0}=
t_{CH}^f(\lambda_1)\wedge t_{CH}^g(\lambda_2)|_{\lambda_2=0}.
\]
Now, by assumption 
\[
(\lambda_1,\lambda_2)\mapsto 
t_\bullet^f(\lambda_1)\wedge t_\bullet^g(\lambda_2),
\]
where $\bullet$ stands for either $CFL$ or $CH$,
is holomorphic at the origin, and thus it follows that
\[
t_{CFL}^f(\lambda)\wedge t_{CFL}^g(\lambda)|_{\lambda=0}=
t_{CH}^f(\lambda)\wedge t_{CH}^g(\lambda)|_{\lambda=0},
\]
but the right hand side is, by \eqref{penna}, equal to
$t_{CFL}^{f\oplus g}(\lambda)|_{\lambda=0}$, and so we obtain Theorem \ref{main_thm} for $r=2$. 
However, the argument easily extends to arbitrary ~$r$.

The way we actually prove Theorem ~\ref{main_thm}, that is, as already announced, by Theorem ~\ref{potential}, is more direct and relies neither on the holomorphicity of the Mellin transform nor on Theorem ~\ref{sondag}. The proof is 
inspired by Proposition ~4.2 in \cite{A}, in which potentials were used to prove Theorem ~\ref{sondag}.
Our hope is that this construction of potentials will be of use for further investigations in the case of a non-complete intersection. 

\begin{proof}[Proof of Theorem ~\ref{potential}]
We let
\begin{equation*}
V =
U^{f_1} +
U^{f_2} \wedge R^{f_1} +
U^{f_3} \wedge R^{f_2} \wedge R^{f_1} +
\ldots +
U^{f_r} \wedge R^{f_{r-1}}\wedge \ldots \wedge R^{f_1}.
\end{equation*}
To motivate this choice of $V$, note that on a formal level
\begin{equation}\label{nablave}
\nabla_f(U^{f_i}\wedge R^{f_{i-1}}\wedge\ldots\wedge R^{f_1})=
R^{f_{i-1}}\wedge\ldots\wedge R^{f_1}-R^{f_{i}}\wedge\ldots\wedge R^{f_1},
\end{equation}
so that
\begin{equation*}\nabla_f V = 
\\ 1-R^{f_{r}}\wedge\ldots\wedge R^{f_1}
.\end{equation*}
Indeed, observe that $\nabla_f$ acts on $U^{f_i}$ just as $\nabla_{f_i}$, so that $\nabla_f U^{f_i}=1-R^{f_i}$. Thus, to prove the first claim of the theorem we have to make this computation legitimate.

First, notice that if a form $A(\lambda)$, depending on a parameter $\lambda$, has an analytic continuation as a current to $\lambda=0$, then clearly $\nabla_f A(\lambda)$ has one.
The action on a test form $\phi$ is given by
\[
\pm\int A(\lambda)\wedge \nabla_f \phi.
\]
However, by integration by parts with respect to $\nabla_f$ and due to the uniqueness of analytic continuations, this is equal to
\[\int\nabla_f A(\lambda)\wedge \phi.\]
To be able to perform the integration by parts in a stringent way we have to regard the currents $T\in\D'_{0,k}(\Lambda^\ell E)$ as functionals on 
\linebreak
$\D_{n,n-k}(\Lambda^{n-\ell}E\wedge \Lambda^n E^*)$. So far we have been a little sloppy about this.

Thus, to compute $\nabla_f V$ we consider the form
\begin{multline*}
v^\lambda = 
|f_1|^{2\lambda} u^{f_1}+
|f_2|^{2\lambda} u^{f_2}\wedge \dbar|f_1|^{2\lambda} \wedge u^{f_1}+ \ldots
\\
\ldots +
|f_r|^{2\lambda} u^{f_r}\wedge \dbar|f_{r-1}|^{2\lambda} \wedge u^{f_{r-1}}
\wedge \ldots \wedge \dbar|f_1|^{2\lambda} \wedge u^{f_1},
\end{multline*}
since, by definition, $v^\lambda|_{\lambda=0}=V$, and accordingly $\nabla_f V=(\nabla_f v^\lambda)|_{\lambda=0}$. More precisely, to verify \eqref{nablave}, let us consider
(recall that $\nabla_f u^{f_i}=1$) 
\begin{multline*}
\nabla_f (|f_i|^{2\lambda} u^{f_i}\wedge \dbar|f_{i-1}|^{2\lambda} \wedge u^{f_{i-1}}\wedge \ldots \wedge \dbar|f_1|^{2\lambda} \wedge u^{f_1}) =\\
-\dbar|f_i|^{2\lambda} \wedge u^{f_i}\wedge \dbar|f_{i-1}|^{2\lambda} \wedge u^{f_{i-1}}\wedge \ldots \wedge \dbar|f_1|^{2\lambda} \wedge u^{f_1}
+\\
|f_i|^{2\lambda}\dbar|f_{i-1}|^{2\lambda} \wedge u^{f_{i-1}}\wedge \ldots \wedge \dbar|f_1|^{2\lambda} \wedge u^{f_1}
+ \mathcal R,
\end{multline*}
where $\mathcal R$ is a sum of terms of the form
\[
|f_i|^{2\lambda} u^{f_i} \wedge 
\dbar|f_{i-1}|^{2\lambda} \wedge u^{f_{i-1}}\wedge \ldots \dbar|f_j|^{2\lambda}
\wedge \dbar|f_{j-1}|^{2\lambda}\wedge \ldots
\wedge \dbar|f_1|^{2\lambda} \wedge u^{f_1},
\]
that arise when $\nabla_f$ falls on any $u^{f_j}$, $j<i$.
The value at $\lambda=0$ of the first term is just 
$-R^{f_i}\wedge R^{f_{i-1}}\wedge\ldots\wedge R^{f_1}$, and it follows from Lemma ~\ref{mymla} that the second term has an analytic continuation to $\lambda=0$ equal to $R^{f_{i-1}}\wedge\ldots\wedge R^{f_1}$. 
The remaining terms, $\mathcal R$, vanish according to Lemma ~\ref{marra}.
Thus \eqref{potential1} is proved, and thereby the first part of the theorem.

Furthermore, let
\begin{multline*}
U^f\wedge V =
U^f\wedge U^{f_1} +
U^f\wedge U^{f_2} \wedge R^{f_1} + \\
U^f\wedge U^{f_3} \wedge R^{f_2} \wedge R^{f_1} +
\ldots +
U^f\wedge U^{f_r} \wedge R^{f_{r-1}}\wedge \ldots \wedge R^{f_1}.
\end{multline*}
We compute $\nabla_f$ of each term. To do this we use a form as above whose analytic continuation to $\lambda=0$ is equal to this particular current. Now, we actually need the extended version of Proposition ~\ref{existensen}, that is Proposition ~\ref{existensen2}. Indeed, consider
\begin{multline*}
\nabla_f(|f|^{2\lambda} u^f \wedge |f_i|^{2\lambda} u^{f_i}\wedge \dbar|f_{i-1}|^{2\lambda} \wedge u^{f_{i-1}}\wedge \ldots \wedge \dbar|f_1|^{2\lambda} \wedge u^{f_1})=
\\
-\dbar|f|^{2\lambda} \wedge u^f \wedge |f_i|^{2\lambda} u^{f_i}\wedge \dbar|f_{i-1}|^{2\lambda} \wedge u^{f_{i-1}}\wedge \ldots \wedge \dbar|f_1|^{2\lambda} \wedge u^{f_1}
+\\
|f|^{2\lambda} |f_i|^{2\lambda} u^{f_i}\wedge \dbar|f_{i-1}|^{2\lambda} \wedge u^{f_{i-1}}\wedge \ldots \wedge \dbar|f_1|^{2\lambda} \wedge u^{f_1}
+\\
|f|^{2\lambda} u^f \wedge \dbar |f_i|^{2\lambda} \wedge u^{f_i}\wedge \dbar|f_{i-1}|^{2\lambda} \wedge u^{f_{i-1}}\wedge \ldots \wedge \dbar|f_1|^{2\lambda} \wedge u^{f_1}
\\
-|f|^{2\lambda} u^f \wedge |f_i|^{2\lambda} \wedge \dbar|f_{i-1}|^{2\lambda} \wedge u^{f_{i-1}}\wedge \ldots \wedge \dbar|f_1|^{2\lambda} \wedge u^{f_1}
+\\
|f|^{2\lambda} u^f \wedge \mathcal R.
\end{multline*}
The first term corresponds to $-R^f\wedge U^{f_i}\wedge R^{f_{i-1}}\wedge \ldots \wedge R^{f_1}$. 
Since $f$ is a complete intersection and $R^f$ therefore is of top degree in $d\bar z_j$ according to Proposition ~\ref{tove}, it is most reasonable to expect also this product to be of top degree in $d\bar z_j$, but because of the factor $U^{f_i}$ $\in \L^{-1}(E_i)$ 
that is apparently not possible unless the product vanishes. 
This is indeed the case, as follows from Lemma ~\ref{snork}. 
The second, third and fourth terms have analytic continuations as 
$U^{f_i}\wedge R^{f_{i-1}}\wedge \ldots \wedge R^{f_1}$,
$U^{f}\wedge R^{f_{i}}\wedge \ldots \wedge R^{f_1}$ and
$-U^{f}\wedge R^{f_{i-1}}\wedge \ldots \wedge R^{f_1}$,
respectively, by Lemma ~\ref{mymla}.
The remaining terms vanish according to Lemma ~\ref{marra}.
Hence 
\begin{multline*}
\nabla_f (U^f\wedge V)=
\sum U^{f_i}\wedge R^{f_{i-1}}\wedge \ldots \wedge R^{f_1} \\-
\sum (U^f\wedge R^{f_{i-1}}\wedge \ldots \wedge R^{f_1} -
U^f\wedge R^{f_{i}}\wedge \ldots \wedge R^{f_1})
=\\
V-U^f+U^f\wedge R^{f_r}\wedge \ldots \wedge R^{f_1}.
\end{multline*}
Finally, the term $U^f\wedge R^{f_r}\wedge \ldots \wedge R^{f_1}$ vanishes by Lemma ~\ref{snus}, and thus taking the lemmas ~\ref{mymla} to ~\ref{snus} for granted, the theorem is proved.
\end{proof}

What remains is the technical part, to prove the lemmas. 
We have tried to put them as simply as possible. Still the formulations may seem a bit strained though. Hopefully, the remarks will shed some light on what matters. We will use the word codegree for the difference between the dimension $n$ of $X$ and the degree.
\begin{lma}\label{mymla}
Let $f=f_1\oplus\cdots\oplus f_r$ be a section of $E^*=E_1^*\oplus\cdots \oplus E_r^*$. Assume that $f$ is a complete intersection. Let $s<r$ and $r'\leq r$. If $h=f$, or if $h=f_i$ for some $i>s$,
then
\begin{equation}\label{form_mymla}
|h|^{2\lambda}|f_{r'}|^{2\lambda}u^{f_{r'}} \wedge 
\ldots \wedge
|f_{s+1}|^{2\lambda}u^{f_{s+1}} \wedge 
\dbar|f_{s}|^{2\lambda}\wedge u^{f_{s}} \wedge 
\ldots \wedge
\dbar|f_1|^{2\lambda}\wedge u^{f_1}
\end{equation}
has an analytic continuation to $\Re\lambda>-\epsilon$, which for $\lambda=0$ is equal to the current $U^{f_{r'}}\wedge\ldots\wedge U^{f_{s+1}}\wedge R^{f_s}\wedge \ldots \wedge R^{f_1}$. 

Moreover,
\begin{equation}\label{form_mymla2}
|h|^{2\lambda}
|f|^{2\lambda}u^{f} \wedge 
|f_{r'}|^{2\lambda}u^{f_{r'}} \wedge 
\ldots \wedge
|f_{s+1}|^{2\lambda}u^{f_{s+1}} \wedge 
\dbar|f_{s}|^{2\lambda}\wedge u^{f_{s}} \wedge 
\ldots \wedge
\dbar|f_1|^{2\lambda}\wedge u^{f_1}
\end{equation}
has an analytic continuation to $\Re\lambda>-\epsilon$, which for $\lambda=0$ is equal to the current $U^f\wedge U^{f_{r'}}\wedge\ldots\wedge U^{f_{s+1}}\wedge R^{f_s}\wedge \ldots \wedge R^{f_1}$. 
\end{lma}

\begin{remark}
The crucial point is that inserting a factor $|h|^{2\lambda}$ has no effect on the value at $\lambda=0$, as long as 
\[
\codim \{h=0\}\cap Y_s\cap\ldots\cap Y_1 > 
\codim Y_s\cap\ldots\cap Y_1,
\]
since then all possibly ``dangerous'' contributions to the current will vanish for degree reasons as in the proof of Proposition ~\ref{tove}. There might be a more general formulation of the lemma that catches this behaviour better. 
Nevertheless, for the proof of Theorem ~\ref{potential} it suffices with the cases above. 
\end{remark}

\begin{proof}
We give a proof of the first claim of the lemma. The second one, concerning \eqref{form_mymla2}, can be proved along the same lines.

For a compactly supported test form $\phi$, we consider
\begin{equation*}
\int
|h|^{2\lambda}
|f_{r'}|^{2\lambda}u^{f_{r'}} \wedge 
\ldots \wedge
|f_{s+1}|^{2\lambda}u^{f_{s+1}} \wedge 
\dbar|f_s|^{2\lambda}\wedge u^{f_s} \wedge 
\ldots \wedge
\dbar|f_1|^{2\lambda}\wedge u^{f_1}
\wedge \phi.
\end{equation*}
After a resolution of singularities as described in the proof of Proposition ~\ref{existensen2}, for $\Re \lambda$ large enough, this integral is equal to a sum of
\begin{multline}\label{int_mymla}
\int
|\mu_h|^{2\lambda}
|\mu_{r'}|^{2\lambda} a_{r'}^{\lambda} 
\frac{\alpha_{r',\ell_{r'}}}{\mu_{r'}^{\ell_{r'}}} \wedge 
\ldots \wedge
|\mu_{s+1}|^{2\lambda} a_{s+1}^{\lambda} 
\frac{\alpha_{s+1,\ell_{s+1}}}{\mu_{s+1}^{\ell_{s+1}}} \wedge
\\
\dbar (|\mu_s|^{2\lambda} a_s^{\lambda}) \wedge 
\frac{\alpha_{s,\ell_s}}{\mu_s^{\ell_s}} \wedge 
\ldots \wedge
\dbar (|\mu_1|^{2\lambda} a_1^{\lambda}) \wedge 
\frac{\alpha_{1,\ell_1}}{\mu_1^{\ell_1}} \wedge
\tilde \phi,
\end{multline}
where the $a_j$'s are strictly positive functions, the $\mu_j$'s are polynomials in some local coordinates $\tau_j$, the $\alpha_{j,\ell_j}$'s are smooth forms and $\tilde\phi$ is as in the proof of Proposition ~\ref{existensen2}.
The existence of the analytic continuation to $\Re \lambda > -\epsilon$ follows from Lemma ~\ref{main_lemma} as before. 

Our aim is to prove that the factor $|h|^{2\lambda}$ does not affect the value at $\lambda=0$.
Let $\sigma$ be one of the coordinate functions $\tau_k$ that divides $\mu_h$.
When expanding each factor $\dbar(|\mu_j|^{2\lambda}a_j^\lambda)$
 by Leibniz' rule we get two different types of terms, integrals with an occurrence of a factor $\dbar|\sigma^\alpha|^{2\lambda}$ for some $\alpha$, and integrals with no such factors.
In the second case the extra factor $|\sigma|^{2\lambda}$ does no harm, since, in fact, the value at $\lambda=0$ is independent of the number of $|\sigma|^{2\lambda}$'s in the numerator as long as there is no $\bar\sigma$ in the denominator. Furthermore, we claim that each integral of the first kind actually vanishes at $\lambda=0$. The argument is analogous to the one in the proof of Proposition ~\ref{tove}. Let us first consider the case when $h=f$.
Observe that the terms in \eqref{form_mymla} are of degree at most $m_1+\ldots + m_{r'} - r'+s \leq m -1$ in $d\bar z_j$, where $m=m_1+\ldots + m_r$. The crucial term -1 appears because of the (at least for the proof) necessary condition that $r<s$, that is that we have at least one factor $U$. Thus, it is enough to consider test forms of 
codegree in $d\bar z$ at most $m-1$. 
We assume that $\phi=\phi_I\wedge d\bar z_I$, where $\phi_I$ is a smooth $(n,0)$-form and $d\bar z_I=d\bar z_{I_1}\wedge \ldots \wedge d\bar z_{I_p}$ where $p\geq n- (m_1+\ldots +m_r) +1$.
Now, $d\bar z_I$ vanishes on the variety $Y=f^{-1}(0)$,
since it has codimension $m$, and accordingly $\Pi^*(d\bar z_I)$ vanishes on $\widetilde Y= \Pi^{-1}Y$, and in particular on $\{\sigma=0\}$. Since it is a form in $d\bar\tau_j$ with antiholomorphic coefficients, each of its terms contains a factor $\bar\sigma$ or $d\bar\sigma$, see the proof of Proposition ~\ref{tove}2, and so in both cases the $\sigma$-integrals vanish according to Lemma ~\ref{main_lemma}.

In the second case, when $h=f_i$, the proof becomes slightly more complicated. We want to prove that the $\sigma$-integral vanishes due to the occurrence of a factor $\bar\sigma$ or $d\bar\sigma$ as above, but
now the desired factors $\bar \sigma$ and $d \bar \sigma$ do not necessarily divide the test form $\tilde\phi$. 
We need to look at a ``larger'' form than ~$\phi$, in fact at the ``largest'' possible ``$\sigma$-free'' form. 
 Without loss of generality we may assume that, for some numbers $s'$ and $r'$, $1 \leq s' \leq s \leq r' \leq r$, $\sigma$ divides $\mu_{s'+1}, \ldots , \mu_s$ and $\mu_{r'+1}, \ldots , \mu_r$ but neither $\mu_1,\ldots ,\mu_s'$ nor $\mu_{s+1}, \ldots, \mu_r'$.
Recall that $u^{f_i}=\sum_\ell v^{f_i}_{\ell_i}/|f_i|^{2\ell}$, where $v^{f_i}_{\ell_i}=s_i\wedge (\dbar s_i)^{\ell_i-1}$.
Let the smooth form
\[
v^{f_r}_{\ell_r}\wedge
\ldots\wedge 
v^{f_{s+1}}_{\ell_{s+1}}\wedge 
\dbar |f_{s'}|^2 \wedge v^{f_s'}_{\ell_s'}\wedge
\ldots\wedge
\dbar |f_{1}|^2 \wedge v^{f_1}_{\ell_1}
\]
be denoted by $F_\ell$, and let 
\[
Y'=\{f_{s'+1}=\ldots =f_s=f_{r'+1}+\ldots f_r=h=0\}.
\]
As above we may assume that $\phi$ consists of only one term $\phi_I\wedge d\bar z_I$. Then, by inspection, the form $F_\ell\wedge d\bar z_I$ is of codegree at most
\[
m_{s'+1}+\ldots +m_s + m_{r'+1}+\ldots + m_r -r + r'
\]
in $d\bar z_j$,
which is strictly less than 
\[
\codim Y'
=m_{s'+1}+\ldots + m_s+m_{r'}+\ldots +m_r+m_h,
\]
because of the assumptions of complete intersection.
Consequently $F_\ell\wedge d\bar z_I$ vanishes
on $Y'$, and thus $\Pi^*(F_\ell \wedge d\bar z_I)$ vanishes on $\Pi^{-1}Y'$, and in particular on $\{\sigma=0\}$. Since it is a form in $d\bar\tau_j$ with antiholomorphic coefficients, each of its terms contains a factor $\bar \sigma$ or a factor $d\bar\sigma$. Using that $\dbar|f|^{2\lambda}=\lambda |f|^{2(\lambda-1)}\dbar |f|^2$, we can write \eqref{int_mymla} as 
\begin{multline*}\label{int_mymla}
\pm \int
|\mu_h|^{2\lambda}
|\mu_r|^{2\lambda} a_r^{\lambda} \frac{\alpha_{r,\ell_r}}{\mu_r^{\ell_r}} \wedge 
\ldots \wedge
|\mu_{r'+1}|^{2\lambda} a_{r'+1}^{\lambda} \frac{\alpha_{r'+1,\ell_{r'+1}}}{\mu_{r'+1}^{\ell_{r'+1}}} \wedge
\\
\dbar (|\mu_s|^{2\lambda} a_s^{\lambda}) \wedge \frac{\alpha_{s,\ell_s}}{\mu_s^{\ell_s}} \wedge 
\ldots \wedge
\dbar (|\mu_{s'+1}|^{2\lambda} a_1^{\lambda}) \wedge \frac{\alpha_{s'+1,\ell_{s'+1}}}{\mu_{s'+1}^{\ell_{s'+1}}} \wedge
\\
\frac
{|\mu_{r'}|^{2\lambda} a_{r'}^{\lambda} \cdots 
|\mu_{1}|^{2\lambda} a_{1}^{\lambda}}
{|\mu_{r'}|^{2\ell_{r'}} \cdots 
|\mu_{s+1}|^{2\ell_{s+1}}}
~~\lambda^{s'}
\frac
{|\mu_{s'}|^{2\lambda} a_{s'}^{(\lambda-1)} \cdots 
|\mu_{1}|^{2\lambda} a_{1}^{(\lambda-1)}}
{|\mu_{s'}|^{2(\ell_{s'}+1)} \cdots 
|\mu_{1}|^{2(\ell_1+1)}}
~~~~ \Pi^*(F_\ell) \wedge \tilde\phi,
\end{multline*} 
where the sign depends on the relation between $r,r',s \text{ and }s'$.
Now the only way a factor $\bar \sigma$ in the numerator (more precisely in $\Pi^*(F_\ell) \wedge \tilde\phi$) could be cancelled out when $\lambda$ is small is by the occurrence of a factor $\bar\sigma$ in one of $\mu_1,\ldots ,\mu_{s'}$, but 
that would obviously contradict the assumption made above.
Hence each term in the integral must contain a factor $\bar\sigma$ or $d\bar\sigma$ independently of the value of $\lambda$ and thus the $\sigma$-integral vanishes according to Lemma ~\ref{main_lemma}.
\end{proof}

\begin{lma}\label{snork}
Let $f=f_1\oplus \cdots \oplus f_r$ be a section of $E^*=E_1^*\oplus \cdots \oplus E_r^*$ of rank $m$ and let $h=f\oplus f'$, where $f'$ is a section of the dual bundle of a holomorphic $m'$-bundle $E'$. Assume that $h$ is a complete intersection. If $r>s$, then
\begin{equation}\label{form_snork}
\dbar|h|^{2\lambda}\wedge u^{h} \wedge
|f_r|^{2\lambda}u^{f_r} \wedge 
\ldots \wedge
|f_{s+1}|^{2\lambda}u^{f_{s+1}} \wedge 
\dbar|f_s|^{2\lambda}\wedge u^{f_s} \wedge 
\ldots \wedge
\dbar |f_1|^{2\lambda}\wedge u^{f_1}
\end{equation}
has an analytic continuation to $\Re\lambda > -\epsilon$ that vanishes at $\lambda=0$.
\end{lma}

\begin{remark}\label{vanish}
Notice that the value at $\lambda=0$ corresponds to the current 
\linebreak
$R^h \wedge U^{f_r}\wedge \ldots \wedge U^{f_{s+1}} \wedge R^{f_s} \wedge \ldots \wedge R^{f_1}$. Since $h$ is a complete intersection, $R^h$ 
is of degree $m+m'$ in $d\bar z_j$
according to Proposition ~\ref{tove}, and therefore it is reasonable to expect also the product to be of degree $m+m'$ in ~$d\bar z_j$. However, since the product contains at least one principal value factor, the degree in $e_j$ must be strictly larger than the degree in $d\bar z_j$, and so, 
the product must vanish. We will see that the assumption that $r>s$ is crucial also for the proof.
\end{remark}

\begin{proof}
After a resolution of singularities as described in the proof of Proposition ~\ref{existensen2}, we can write \eqref{form_snork} integrated against a test form $\phi$ as a sum of terms of the type
\begin{multline}\label{int_snork}
\int
\dbar(|\mu_h|^{2\lambda} a_h^{\lambda}) \frac{\alpha_{h,\ell_h}}{\mu_h^{\ell_h}} \wedge 
|\mu_r|^{2\lambda} a_r^{\lambda} \frac{\alpha_{r,\ell_r}}{\mu_r^{\ell_r}} \wedge 
\ldots \wedge
|\mu_{s+1}|^{2\lambda} a_{s+1}^{\lambda} \frac{\alpha_{s+1,\ell_{s+1}}}{\mu_{s+1}^{\ell_{s+1}}} \wedge
\\
\dbar (|\mu_s|^{2\lambda} a_s^{\lambda}) \wedge \frac{\alpha_{s,\ell_s}}{\mu_s^{\ell_s}} \wedge 
\ldots \wedge
\dbar(|\mu_1|^{2\lambda} a_1^{\lambda}) \wedge \frac{\alpha_{1,\ell_1}}{\mu_1^{\ell_1}} \wedge
\tilde \phi,
\end{multline}
where the $\alpha_{i,\ell_i}$'s are smooth forms of degree $\ell_i$ in $e_j$, 
the $a_i$'s are non-vanishing functions and the $\mu_i$'s are monomials in some local coordinates ~$\tau_k$ and $\tilde\phi$ is as in the previous proofs.

We expand the factor $\dbar(|\mu_h|^{2\lambda}a_h^\lambda)$ by Leibniz' rule and consider the term obtained when $\dbar$ falls on $|\sigma|^{2\lambda}$, where $\sigma$ is one of the $\tau_k$'s that divide $\mu_h$. We prove that this term vanishes when integrating with respect to $\sigma$. The term that arises when $\dbar$ falls on $a_h^\lambda$ clearly vanishes as before, see the proof of Proposition ~\ref{existensen2}.
Since the rank of $E\oplus E'$ is $m+m'$, the terms in \eqref{form_snork} are of degree at most $m+m'-1$ in $d\bar z$,
since we have at least one 
$U$-factor. Thus it is enough to consider test forms of codegree in $d\bar z$ at most $m+m'-1$. As in the previous proofs we may assume that $\phi=\phi_I\wedge d\bar z_I$. It follows that ~$ d\bar z_I$ 
\linebreak
vanishes on $Y=h^{-1}(0)$ for degree reasons, and thus ~$\Pi^*(d\bar z_I)$ vanishes on ~$\Pi^{-1} Y$. Since this is a form in $d\bar\tau_j$ with antiholomorphic coefficients, each of its terms contains a factor $\bar \sigma$ or $d\bar\sigma$
and consequently 
the $\sigma$-integral vanishes according to Lemma ~\ref{main_lemma}.
\end{proof}

\begin{lma}\label{marra}
Let $f=f_1\oplus\cdots\oplus f_r$ be a section of $E^*=E_1^*\oplus\cdots\oplus E_r^*$. Assume that $f$ is a complete intersection and let $s<r$. Then
\begin{equation}\label{form_marra}
|f_r|^{2\lambda}u^{f_r} \wedge 
\ldots \wedge
|f_{s+1}|^{2\lambda}u^{f_{s+1}} \wedge 
\dbar|f_s|^{2\lambda}\wedge u^{f_s} \wedge 
\ldots \wedge
\dbar|f_t|^{2\lambda}\wedge
\ldots \wedge
\dbar|f_1|^{2\lambda}\wedge u^{f_1}
\end{equation}
and
\begin{multline}\label{form_marra2}
|f|^{2\lambda}u^f\wedge
|f_r|^{2\lambda}u^{f_r} \wedge 
\ldots \wedge
|f_{s+1}|^{2\lambda}u^{f_{s+1}} \wedge 
\dbar|f_s|^{2\lambda}\wedge u^{f_s} \wedge \\
\ldots \wedge
\dbar|f_t|^{2\lambda}\wedge
\ldots \wedge
\dbar|f_1|^{2\lambda}\wedge u^{f_1}
\end{multline}
have analytic continuations to $\Re\lambda>-\epsilon$ that vanish at $\lambda=0$.
\end{lma}

\begin{remark}
Morally, what this lemma says is that when applying Leibniz' rule to $\nabla_f$ acting on a product of principal value and residue currents, there will be no contributions from $\nabla_f$ falling on a residue factor. Of course this is expected, since the residue currents are $\nabla_f$-closed.
\end{remark}

\begin{proof}
For \eqref{form_marra} the result follows from Lemma ~\ref{mymla} after an integration by parts with respect to $\nabla_f$. (Recall that $\phi$ is a form taking values in $\Lambda^{n-\ell}E\wedge \Lambda^n E^*$.) Note that $\dbar|f_t|^{2\lambda}=-\nabla_f |f_t|^{2\lambda}$.
By Stokes' theorem,
\begin{multline*}
\int
|f_r|^{2\lambda}u^{f_r} \wedge 
\ldots \wedge
|f_{s+1}|^{2\lambda}u^{f_{s+1}} \wedge 
\dbar|f_s|^{2\lambda}\wedge u^{f_s} \wedge \ldots \\
\ldots \wedge
\nabla_f|f_t|^{2\lambda}\wedge\ldots
\wedge
\dbar|f_1|^{2\lambda}\wedge u^{f_1}
\wedge \phi
=\\
\pm\int
(|f_t|^{2\lambda}-1)
\nabla_f(|f_r|^{2\lambda}u^{f_r} \wedge 
\ldots \wedge
|f_{s+1}|^{2\lambda}u^{f_{s+1}} \wedge 
\\
\dbar|f_s|^{2\lambda}\wedge u^{f_s} \wedge 
\ldots \wedge
\dbar|f_1|^{2\lambda}\wedge u^{f_1}
\wedge \phi),
\end{multline*} 
so it is enough to prove that this expression vanishes at $\lambda=0$.
For degree reasons it will vanish if $\nabla_f$ acts on $\phi$, so let us consider the form
\begin{equation}\label{rosie}
\nabla_f(|f_r|^{2\lambda}u^{f_r} \wedge 
\ldots \wedge
|f_{s+1}|^{2\lambda}u^{f_{s+1}} \wedge 
\dbar|f_s|^{2\lambda}\wedge u^{f_s} \wedge 
\ldots \wedge
\dbar|f_1|^{2\lambda}\wedge u^{f_1}).
\end{equation}
Now, applying Leibniz' rule gives a sum of terms, of which the ones arising when $\nabla_f$ falls on a principal value term will vanish for degree reasons, whereas the others will be precisely as in the hypothesis of Lemma ~\ref{mymla}. Moreover $f_t$ is an $h$ of the second kind, so according to Lemma ~\ref{mymla} the factor $|f_t|^{2\lambda}$ does not have any effect on the value at $\lambda=0$. Thus we are done.

In the case of \eqref{form_marra2}, after an integration by parts, we have to prove that
\begin{multline*}
\int
(|f_t|^{2\lambda}-1)
\nabla_f(|f|^{2\lambda}u^{f} \wedge
|f_r|^{2\lambda}u^{f_r} \wedge 
\ldots \wedge
|f_{s+1}|^{2\lambda}u^{f_{s+1}} \wedge 
\\
\dbar|f_s|^{2\lambda}\wedge u^{f_s} \wedge 
\ldots \wedge
\dbar|f_1|^{2\lambda}\wedge u^{f_1}
\wedge \phi)
\end{multline*} 
vanishes at $\lambda=0$.
The term when $\nabla_f$ falls on the factor $|f|^{2\lambda}u^f$ is of the type in Lemma ~\ref{snork}. It is easy to see from the proof that the factor ~$|f_t|^{2\lambda}$ does not affect the value at $\lambda=0$ and so this term vanishes. The remaining part is as in the hypothesis of the latter statement of Lemma ~\ref{mymla}, thus the result follows as above. 
\end{proof}

\begin{lma}\label{snus}
Let $f=f_1\oplus\cdots\oplus f_r$ be a section of $E^*=E_1^*\oplus \cdots \oplus E_r^*$. Assume that $f$ is a complete intersection. Let $h=f_{I_1}\oplus \cdots \oplus f_{I_p}$, where $I=\{I_1,\ldots, I_p \} \subset \{1,\ldots, r\}$.
Then
\begin{equation}\label{form_snus}
|h|^{2\lambda}u^{h} \wedge
\dbar|f_r|^{2\lambda}\wedge u^{f_r} \wedge 
\ldots \wedge
\dbar |f_1|^{2\lambda}\wedge u^{f_1}
\end{equation}
has an analytic continuation to $\Re\lambda > -\epsilon$ that vanishes at $\lambda=0$.
\end{lma}

\begin{remark}
The value at $\lambda=0$ corresponds to the current $U^h\wedge R^{f_1}\wedge \ldots\wedge R^{f_r}$. Since the $R$-part is of top degree according to Proposition ~\ref{tove} this product should formally vanish by arguments similar to those in Remark ~\ref{vanish}.
\end{remark}

\begin{proof}
As in the proofs of the previous lemmas we start by a resolution of singularities.
Thus, the form 
\eqref{form_snus} integrated against a test form $\phi$ is equal to a sum of terms of the type
\begin{multline}\label{int_snus}
\int
|\mu_h|^{2\lambda} a_h^{\lambda} \frac{\alpha_{h,\ell_h}}{\mu_h^{\ell_h}} \wedge 
\dbar(|\mu_r|^{2\lambda} a_r^{\lambda}) \wedge \frac{\alpha_{r,\ell_r}}{\mu_r^{\ell_r}} \wedge 
\ldots \wedge
\dbar(|\mu_1|^{2\lambda} a_1^{\lambda}) \wedge \frac{\alpha_{1,\ell_1}}{\mu_1^{\ell_1}} \wedge
\tilde \phi,
\end{multline}
where $\alpha_{i,\ell_i}$, $a_i$, $\mu_i$ and $\tilde\phi$ are as above.
Further, we can find a resolution to a certain toric variety so that locally one of the monomials $\mu_1, \ldots ,\mu_r$ divides the other ones. Without loss of generality we may assume that ~$\mu_1$ divides 
all $\mu_j$'s.
We expand $\dbar(|\mu_1|^{2\lambda}a_1^\lambda)$ by Leibniz' rule. The term obtained when ~$\dbar$ falls on $a_1^\lambda$ vanishes as in the proof of Proposition ~\ref{existensen2}, so it is enough to consider the terms that arise when ~$\dbar$ falls on $|\sigma|^{2\lambda}$, where $\sigma$ is one of the coordinates in $\mu_1$.

We claim that the $\sigma$-integral vanishes at $\lambda=0$. As usual, we observe that the terms of \eqref{form_snus} are of degree at most $m-1$ in $d\bar z_j$, where the -1 in this case is due to the factor $U^h$, so it suffices to consider test forms of codegree at most $m-1$. We assume that $\phi=\phi_I\wedge d\bar z_I$, where $\phi_I$ is an $(n,0)$-form and $d\bar z_I=d\bar z_{I_1} \wedge \ldots \wedge d\bar z_{I_p}$, where $p\leq n-m+1$. Then $d\bar z_I$ vanishes on 
the variety $Y=f^{-1}(0)$ for degree reasons, and accordingly $\Pi^*(d\bar z)$ vanishes on $\Pi^{-1}Y$, and in particular on $\{\sigma_1=0\}$. By arguments as in the proof of Proposition ~\ref{tove} it follows that $\Pi^*(d\bar z)$ must contain a factor $\bar\sigma$ or $d\bar\sigma$ since it is a form in $d\bar \tau_k$ with antiholomorphic coefficients, and hence the $\sigma$-integral 
vanishes as before.
\end{proof}

\section{An example}\label{examples}
We conclude this paper with an explicit computation, by which we enlighten the possibility of extending Theorem ~\ref{main_thm} to a slightly weaker notion of complete intersection.
Indeed, when generalizing Theorem ~\ref{sondag}, or rather its line bundle formulation ~\eqref{sondag-bundle}, to sections of bundles of arbitrary rank, it is not obvious how one should interpret the assumption of $f$ being a complete intersection. In the formulation of Theorem ~\ref{main_thm} we require 
the codimension of $f^{-1}(0)$ to be equal to the rank of the bundle $E$. A less strong hypothesis would be to just demand the $f_i$'s to intersect properly, that is, that $\codim f^{-1}(0)=p_1+\ldots +p_r$ if $p_i=\codim f_i$. However, the following example shows that Theorem ~\ref{main_thm} does not extend to this case.

\begin{ex}\label{noncomplete}
Let $f_1=z_1^2$, $f_2=z_1 z_2$ and $g=z_2 z_3$. Then 
\begin{multline*}
Y_f=f^{-1}(0)=\{z_1=0\},\,\,\,\,\,\,\,\, Y_g=g^{-1}(0)=\{z_2=0\}\cup\{z_3=0\},
\\\text{ and }
Y=Y_f\cap Y_g=\{z_1=z_2=0\}\cup\{z_1=z_3=0\}.
\end{multline*}
Note that $Y_f$ and $Y_g$ have codimension 1, and that $Y$ has codimension ~2. Thus $f$ and $g$ intersect properly, although they do not define a complete intersection. 

Let us compute $(R^f\wedge R^g)_2$. Adopting the trivial metric we get
\[
s^f=\bar f_1 e_1 + \bar f_2 e_2= \bar z_1( \bar z_1 e_1 + \bar z_2 e_2)
\text{ and }
|f|^2=|z_1|^2(|z_1|^2+|z_2|^2),
\]
so that
\[
u^f_1= \frac{ \bar z_1 e_1 + \bar z_2 e_2}
{z_1(|z_1|^2+|z_2|^2)}.
\]
Let $\phi$ be a test form of bidegree $(3,1)$ with support outside 
$\{z_2=0\}$.
Then $R^f\wedge R^g.\phi$ is given by
\begin{multline*}
\int
\dbar(|z_1|^{2\lambda}(|z_1|^2+|z_2|^2)^\lambda)\wedge
\frac{ \bar z_1 e_1 + \bar z_2 e_2}
{z_1(|z_1|^2+|z_2|^2)}\wedge
\dbar|z_2 z_3|^{2\lambda}\wedge\frac{e_3}{z_2 z_3}\wedge \phi|_{\lambda=0} =\\
-\int
\dbar \Big[\frac{1}{z_1}\Big]
\Big [\frac{1}{z_2}\Big ]\wedge
\dbar\Big[\frac{1}{z_3}\Big]\wedge e_2 \wedge e_3\wedge \phi.
\end{multline*}
Note that the support of the current is on the $z_2$-axis, as expected since ~$\phi$ has support outside the $z_3$-axis. 

To deal with test forms with support intersecting $\{z_2=0\}$ we need to resolve the singularity of $f$ at the $z_3$-axis. Let $\widetilde\U$ be the blow-up of ~$\mathbb C^3_z$ along the $z_3$-axis and let $\Pi:\widetilde\U\to\mathbb C^3_z$ be the corresponding proper map. 
We can cover $\widetilde\U$ by two coordinate charts,
\begin{multline*}
\Omega_1=\{(\tau_1,\tau_2,z_3);\ (\tau_1,\tau_1\tau_2,z_3)=z\in\mathbb C^3_z\}\ \\
\text{and }
\Omega_2=\{(\sigma_1,\sigma_2,z_3);\ (\sigma_1\sigma_2,\sigma_1,z_3)=z\in\mathbb C^3_z\}.
\end{multline*}
In $\Omega_1$
\[
\Pi^* u^f_1=\frac{e_1+\bar \tau_2 e_2}{\tau_1^2(1+|\tau_2|^2)},
\]
and thus
\begin{multline}\label{syrsa}
R^{\Pi^* f}\wedge R^{\Pi^* g}=
\dbar |\tau_1|^{4\lambda} \wedge \frac{e_1+\bar \tau_2 e_2}{\tau_1^2(1+|\tau_2|^2)} \wedge 
~~\dbar |\tau_1 \tau_2 z_3|^{2\lambda} \wedge \frac {e_3}{\tau_1 \tau_2 z_3} \Big |_{\lambda=0}=\\
\frac{2}{3}~
\dbar \Big [\frac{1}{\tau_1^3}\Big] \wedge \frac{e_1+\bar \tau_2 e_2}{\tau_1^3(1+|\tau_2|^2)} \wedge 
\dbar \Big [\frac{1}{\tau_2z_3}\Big ] \wedge e_3.
\end{multline}
Let $\phi$ be a test form of bidegree $(3,1)$; we can write $\phi$ as $\phi_I \wedge dz$, where $dz= dz_1\wedge dz_2 \wedge dz_3$ and 
$\phi_I=\varphi^1(z) d\bar z_1 + \varphi^2(z) d\bar z_2 +\varphi^3(z) d\bar z_3$.
Now $\int R^f\wedge R^g\wedge\phi=\int_{\widetilde \U} R^{\Pi^*f}\wedge R^{\Pi^*g}\wedge\Pi^*\phi$. To compute the contribution from the chart $\Omega_1$, let
$\widetilde\phi=\chi\Pi^* \phi$, where $\chi$ is a function of some partition of unity with support in $\Omega_1$.  We may without loss of generality assume that $\chi$ only depends on $|\tau_i|$ and also that $\chi(0,0,z_3)=1$. 
Then we get that
$R^{\Pi^* f}\wedge R^{\Pi^* g}.\widetilde\phi$ is equal to
\begin{multline*}
\frac{2}{3}\int
\dbar \Big [\frac{1}{\tau_1^3}\Big] \wedge e_1 \wedge 
\dbar \Big [\frac{1}{\tau_2}\Big ] \Big [\frac {1}{z_3} \Big ] \wedge e_3 
\wedge \chi \varphi^3(\tau_1,\tau_1 \tau_2, z_3) ~~d\bar z_3
\wedge d\tau_1 \wedge d(\tau_1\tau_2) \wedge dz_3
=\\
(2\pi i )^2 \frac{2}{3}\int
\varphi^3_{1}(0,0,z_3)\frac{1}{z_3} e_1\wedge e_3\wedge d\bar z_3\wedge dz_3
=\\
- \frac{2}{3}\int
\dbar \Big [\frac{1}{z_1^2}\Big] \wedge 
\dbar \Big [\frac{1}{z_2}\Big ] \Big [\frac {1}{z_3} \Big ] \wedge e_1 \wedge e_3 
\wedge  \phi,
\end{multline*}
where we have used the well known fact that
\begin{equation}\label{salt}
\int_z
\dbar\Big[\frac{1}{z^p}\Big] \wedge \psi(z) dz = 
\frac{2\pi i}{(p-1)!}\frac{\partial^{p-1}}{\partial z^{p-1}}\psi(0).
\end{equation}

Computing $R^{\Pi^*f}\wedge R^{\Pi^*g}$ in $\Omega_2$ gives yet another contribution. Altogether we get
\begin{multline*}
(R^f\wedge R^g)_2=
-\dbar \Big[\frac{1}{z_1}\Big]
\Big [\frac{1}{z_2}\Big ]\wedge
\dbar\Big[\frac{1}{z_3}\Big]\wedge e_2 \wedge e_3\\
- \frac{2}{3}
\dbar \Big [\frac{1}{z_1^2}\Big] \wedge 
\dbar \Big [\frac{1}{z_2}\Big ] \Big [\frac {1}{z_3} \Big ] \wedge e_1 \wedge e_3 +
\frac{1}{3}
\dbar \Big [\frac{1}{z_1}\Big] \wedge 
\dbar \Big [\frac{1}{z_2^2}\Big ] \Big [\frac {1}{z_3} \Big ] \wedge e_2 \wedge e_3.
\end{multline*}
Similar computations yield
\begin{multline*}
R^{f\oplus g}_2
=
\dbar \Big [\frac{1}{z_1^2} \Big ] \wedge \dbar\Big [\frac{1}{z_2} \Big ]
\Big [\frac{1}{z_3} \Big ] \wedge e_1 \wedge e_3
- \dbar \Big [\frac{1}{z_1} \Big ] \Big [\frac{1}{z_2^2} \Big ]\wedge 
\dbar \Big [\frac{1}{z_3} \Big ] \wedge e_2 \wedge e_3.
\end{multline*}
For details, we refer to \cite{W}. See also Theorem 5.2 in \cite{W2}.

To conclude, $R^f\wedge R^g\neq R^{f\oplus g}$, and hence Theorem ~\ref{main_thm} does not generalize to the case of proper intersections.
\end{ex}

{\bf Acknowledgement:} The author would like to thank Mats Andersson, who came up with the original idea of this paper, for many valuable discussions on the topic and helpful comments on preliminary versions.

\def\listing#1#2#3{{\sc #1}:\ {\it #2},\ #3.}

\end{document}